\documentclass[11pt]{article}
\usepackage{latexsym}
\usepackage{amssymb,latexsym, mathptmx}
\usepackage{amsmath}
\usepackage{mathrsfs}
\usepackage{amscd}
\usepackage{amsfonts}
\DeclareFontFamily{OT1}{pzc}{}
\DeclareFontShape{OT1}{pzc}{m}{it}{<-> s*[1.11] pzcmi7t}{}
\DeclareMathAlphabet{\mathcal}{OT1}{pzc}{m}{it}

\newcommand{\ra}{\rightarrow}

\newcommand{\R}{\mathbb{R}}
\newcommand{\di}{\mathrm{d}}

\newcommand{\rad}{\mathrm{rad}}
\newcommand{\grad}{{{\rm{grad}\,}}}
\newcommand{\HessM}{{\rm \Hess}_{\!_M}}

\newcommand{\Hess}{{{\rm Hess}\,}}

\newtheorem{theorem}{\bf Theorem}[section]
\newtheorem{lemma}[theorem]{\bf Lemma}
\newtheorem{corollary}[theorem]{\bf Corollary}

\newtheorem{remark}[theorem]{\bf Remark}
\newtheorem{definition}[theorem]{\bf Definition}

\begin{document}
\hyphenpenalty=300
\title{Curvature estimates for properly immersed $\phi_{h}$-bounded submanifolds}
\author{G. Pacelli Bessa\thanks{The authors were partially supported by   CNPq-Brazil, PROCAD-PICME-CAPES-Brazil and The Abdus Salam International Centre for Theoretical Physics - ICTP, Italy.} \and   Barnabe P. Lima \and  Leandro F. Pessoa}
\date{}
\maketitle

\begin{abstract}Jorge-Koutrofiotis \cite{jorge-koutrofiotis} \& Pigola-Rigoli-Setti \cite{prs-memoirs} proved sharp sectional curvature estimates for extrinsically bounded submanifolds.  Alias,  Bessa and Montenegro in \cite{alias-bessa-montenegro}, showed that these estimates  hold on  properly immersed cylindrically bounded submanifolds. On the other hand, in \cite{alias-bessa-dajczer}, Alias, Bessa and Dajczer  proved sharp mean curvature estimates for properly immersed cylindrically bounded submanifolds.  In this paper we prove these sectional and mean curvature estimates for  a larger class of submanifolds, the properly immersed $\phi$-bounded submanifolds, Thms. \ref{thmMain-intro} \& \ref{thmMain-2}. Thse ideas, in fact,  we  prove stronger forms  of these estimates, see the results in section \ref{sec:OM-Pairs}.

\vspace{2mm}

\noindent \textbf{keywords:}  Curvature estimates,  $\phi$-bounded submanifolds,  Omori-Yau pairs, Omori-Yau maximum principle.

\end{abstract}

\section{Introduction}
\label{intro}The classical isometric immersion problem asks whether there exists an isometric immersion
$\varphi \colon M\to  N$ for given Riemannian manifolds $M$ and $N$ of dimension $m$ and $n$ respectively, with $m<n$. The model result for this type of problem is the celebrated  Efimov-Hilbert Theorem  \cite{efimov}, \cite{hilbert} that says that there is no isometric immersion of a geodesically complete surface $M$ with sectional curvature $K_{M}\leq -\delta^2<0$ into $\mathbb{R}^{3}$,  $\delta\in \mathbb{R}$.  On the other hand, the Nash Embedding Theorem shows  that there is always an isometric embedding into the  Euclidean $n$-space $\mathbb{R}^{n}$  provided  the codimension $n-m$ is sufficiently large, see \cite{nash}.
For small codimension, meaning in this paper that $n-m\leq m-1 $,  the answer in general depends on the geometries of $M$ and $N$. For instance, a classical   result of
C. Tompkins \cite{tompkins} states that a compact, flat, $m$-dimensional Riemannian manifold can not be isometrically immersed into $\mathbb{R}^{2m-1}$. C. Tompkin's  result was extended in  a series of papers, by  Chern and Kuiper
\cite{chern-kuiper}, Moore \cite{moore}, O'Neill \cite{oneil},  Otsuki \cite{otsuki} and Stiel \cite{stiel},  whose results  can
be summarized in the following theorem.

\begin{theorem}[C. Tompkins et al.]\label{thmTompkins}Let $\varphi\colon M\to N$ be an isometric immersion of compact Riemannian $m$-manifold $M$ into  a Cartan-Hadamard $n$-manifold $N$ with small codimension $n-m\leq m-1$.
Then the sectional curvatures of $M$ and $N$ satisfy \begin{equation}\sup_{M} K_{M}> \inf_{N} K_{N}. \label{eqTompkins}\end{equation}
\end{theorem}

  L. Jorge and   D. Koutrofiotis  \cite{jorge-koutrofiotis},  considered complete extrinsically bounded\footnote{Meaning: immersed into regular geodesic balls of a Riemannian manifold.} submanifolds with  scalar curvature bounded from below  and proved the curvature estimates \eqref{eqJorge-koutrofiotis}.
    Pigola, Rigoli and Setti \cite{prs-memoirs} proved an all general and abstract version of the Omori-Yau maximum principle \cite{cheng-yau}, \cite{yau} and in consequence they were able to extend Jorge-Koutrofiotis' Theorem  to  complete  $m$-submanifolds $M$ immersed into  regular balls  of any Riemannian $n$-manifold $N$  with scalar curvature  bounded below  as $s_{_M}\geq - c\cdot\rho^{2}_{_M}\cdot\prod_{j=1}^{k}\Big(\log^{(j)}(\rho_{_M})\Big)^2,\,\, \rho_{_M}\gg 1$.

    Their version of Jorge-Koutrofiotis Theorem is the following.

\begin{theorem}[Jorge-Koutrofiotis \& Pigola-Rigoli-Setti]\label{thmJK}Let $\varphi\colon M\to N$ be an isometric immersion of a complete Riemannian  $m$-manifold $M$ into a $n$-manifold $N$, with $n-m \leq m-1$,  with $\varphi(M)\subset B_N(r)$, where $B_{N}(r)$
is a regular  geodesic ball of $N$.
If the scalar curvature of $M$ satisfies \begin{equation}
s_{_M}\geq -c\cdot \rho^{2}_{_M}\cdot\prod_{j=1}^{k}\Big(\log^{(j)}(\rho_{_M})\Big)^2,\,\, \rho_{_M}\gg 1,
\label{eqScalar}
\end{equation}
for some constant $c>0$ and some integer $k\geq 1$,
where $\rho_{_M}$ is the distance function on $M$ to a fixed point and $\log^{(j)}$ is the $j$-th iterate of the
logarithm. Then
\begin{equation}\label{eqJorge-koutrofiotis}
\sup_{M}K_{M}\geq C_{b}^{2}(r)+\inf_{B_{N}(r)}K_{N},
\end{equation}
where $b=\sup_{B_{_N}(r)}K_{N}^{\rad}\leq b$
\begin{equation}\label{eqCb}
C_b(t)=
\begin{cases}
\sqrt{b}\cot(\sqrt{b}\,t) & \mbox{\rm if $b>0$ and $0<t<\pi/2\sqrt{b}$}\\
1/t & \mbox{\rm if $b=0$ and $t>0$}\\
\sqrt{-b}\coth(\sqrt{-b}\,t) & \mbox{\rm if $b<0$ and $t>0$}.
\end{cases}
\end{equation}
\end{theorem}
\begin{remark}If $B(r)\subset \mathbb{N}^{n}(b)$ is a geodesic ball of radius $r$ in  the simply connected space form    of sectional curvature $b$,
$\partial B(r)$  its boundary and $\varphi \colon \partial B(r-\epsilon) \to B(r)$ is the canonical immersion, where $\epsilon>0$ is small, then we have  \[
\sup_M K_M=K_{\partial B(r-\epsilon)}=
\begin{cases}
{b}/{\sin^2(\sqrt{b}\,(r-\epsilon))} & \mbox{\rm if $b>0$} \\
1/(r-\epsilon)^2 & \mbox{\rm if $b=0$ }\\
{-b}/\sinh^2(\sqrt{-b}\,(r-\epsilon)) & \mbox{\rm if $b<0$ }.
\end{cases}
\]
Therefore,
$
\sup_M K_M-[C_{b}^{2}(r)+\inf K_{\mathbb{N}^{n}(b)}]= [C_b^2(r-\epsilon)-C_b^2(r)]\to 0$ as $\epsilon \to 0$,  showing that   the inequality
(\ref{eqJorge-koutrofiotis}) is sharp.
\end{remark}
\begin{remark}One may assume without loss of generality that $\sup_M K_M<\infty$. This together with the scalar curvature bounds \eqref{eqScalar} implies  that $$K_{M}\geq - c^{2}\cdot \rho^{2}_{_M}\cdot\prod_{j=1}^{k}\Big(\log^{(j)}(\rho_{_M})\Big)^2,\,\, \rho_{_M}\gg 1$$ for some positive constant $c>0$. This curvature lower bound implies that $M$ is stochastically complete,  which  it is equivalent  to the fact that $M$ hold the weak maximum principle, (a weaker form of Omori-Yau maximum principle, see details in \cite{PRS-PAMS}),  and that is enough to reproduce Jorge-Koutrofitis original proof of the curvature estimate  \eqref{eqJorge-koutrofiotis}.  \label{remark4}
\end{remark}

Recently,  Alias,  Bessa and Montenegro   \cite{alias-bessa-montenegro}  extended  Theorem \ref{thmJK}  to the class of cylindrically bounded, properly immersed submanifolds, where an isometric immersion $\varphi \colon M\hookrightarrow N\times \mathbb{R}^{\ell}$ is   cylindrically bounded if $\varphi (M)\subset B_{N}(r)\times \mathbb{R}^{\ell}$.  Here $B_{N}(r)$ is a geodesic ball in $N$ of radius $r>0$. They proved the following theorem.

\begin{theorem}[Alias-Bessa-Montenegro]\label{thmABM} Let $\varphi \colon M\rightarrow N\times\R^{\ell}$ be a  cylindrically bounded isometric immersion,  $\varphi(M)\subset B_N(r)\times\R^{\ell}$, where $B_N(r)$ is a regular geodesic ball of $N$ and $b=\sup K^{\rad}_{B_{N}(r)}$.   Let ${\rm dim}(M)=m$, ${\rm dim}(N)=n-\ell$ and assume that $n-m\leq m-\ell -1$.  If either
\begin{itemize}
\item[i.] the scalar curvature of $M$ is bounded below as  (\ref{eqScalar}), or
\item[]
\item[ii.] the immersion $\varphi$ is proper and
\begin{equation}
\sup_{\varphi^{-1}(B_{N}(r)\times \partial B_{\mathbb{R}^{\ell}}(t))}\Vert\alpha\Vert\leq \sigma(t),
\label{growth}
\end{equation}
\end{itemize}
where $\alpha$ is the second fundamental form of $\varphi$ and $\sigma:[0,+\infty)\rightarrow\mathbb{R}$ is a
positive  function satisfying $\int_0^{+\infty}dt/\sigma(t)=+\infty$, then
\begin{equation}
\sup_{M}K_{M}\geq C_{b}^{2}(r)+\inf_{B_{N}(r)}K_{N}.
\label{eq-ABM}
\end{equation}
\end{theorem}
\begin{remark} The idea  is to show that  the hypotheses, in both items i. $\&$ ii.  implies that $M$ is stochastically complete, then  Remark \ref{remark4} applies.
\end{remark}
In the same spirit,
 Alias, Bessa and Dajczer \cite{alias-bessa-dajczer}, had proved the following mean curvature estimates for  cylindrically bounded submanifolds properly immersed into $N\times \mathbb{R}^{\ell}$ immersed submanifolds.
\begin{theorem}[Alias-Bessa-Dajczer]\label{thm-Alias-Bessa-Dajczer}Let $\varphi \colon M\rightarrow N\times\R^{\ell}$ be a  cylindrically bounded isometric immersion,  $\varphi(M)\subset B_N(r)\times\R^{\ell}$, where $B_N(r)$ is a regular geodesic ball of $N$ and $b=\sup K^{\rad}_{B_{N}(r)}$. Here
$M$ and $N$ are complete Riemannian manifolds of dimension $m$ and $n-\ell$ respectively, satisfying $m\geq \ell+1$.
 If the immersion $\varphi$ is proper, then
\begin{equation}\label{eqABD-mean}
\sup_{M}\vert H\vert \geq (m-\ell)\cdot C_{b}(r).
\end{equation}
\end{theorem}
\section{ Main results}
The purpose of this paper  is to extend these curvature estimates to a larger class  of submanifolds, precisely, the properly immersed $\phi$-bounded submanifolds. To describe this class we  need to  introduce  few preliminaries.
\subsection{ $\phi$-bounded submanifolds}\label{sec:phi-bounded}
Consider
 $G \in C^{\infty}([0, \infty))$  satisfying \begin{equation}\label{buonacurva2}G_{-}\in L^{1}(\mathbb{R}^{+}),\,\,\,\,\,
t \int_{t}^{+\infty} G_-(s) \di s \le \frac{1}{4} \,\,\,\, \text{on }\,\, \R^+,
\end{equation} and  $h$  the solution of the following differential equation
\begin{align}\label{eqg}
\left\{
  \begin{array}{l}
   h''(t)-G(t)h(t)=0, \\[0.1cm]
   h(0)=0, \,\,\,\,\; h'(0)=1.
  \end{array}
\right.
\end{align} In \cite[Prop.  1.21]{bianchini-mari-rigoli}, it is proved  that the solution $h$  and its derivative $h'$ are positive  in $\R^+=(0, \infty)$, provided $G$ satisfies \eqref{buonacurva2}
and furthermore $h \ra +\infty$ whenever the stronger condition
\begin{equation}\label{buonacurva}
G(s) \ge -\frac{1}{4s^2} \qquad \text{on } \R^+
\end{equation}
holds. Define $\phi_h\in C^{\infty}([0, \infty))$  by \begin{equation}\phi_{h}(t)=\int_{0}^{t}h(s)ds.\label{eq-phi}\end{equation}Since $h$ is positive and increasing in $\R^+$,   we have that   $\lim_{t\to \infty}\phi_h(t)=+\infty$. Moreover, $\phi_h$ satisfies the differential equation $$\phi_h''(t) - \displaystyle\frac{h'}{h}(t)\phi_h'(t)=0$$ for all $t\in [0, \infty)$.
\vspace{2mm}
\paragraph{Notation.} In  this paper, $N$ will always be  a complete  Riemannian manifold with a  distinguished point $z_0$ and radial sectional curvatures along the minimal geodesic  issuing from $z_0$ bounded above by $K_{N}^{rad} (z)\leq -G(\rho_N(z))$, where  $G$ satisfies the conditions \eqref{buonacurva2}. Let    $h$ be the solution of \eqref{eqg} associated to $G$  and   $\phi_{h}=\int h(s)ds$. Finally,  $\rho_{N}(z)={\rm dist}_{N}(z_0,z)$ will be the distance function on $N$. For any given    complete Riemannian manifold $(L,y_0)$ with a distinguished point $y_0$ and radial sectional curvature\footnote{Along the geodesics issuing from $y_0$.} bounded below ($K^{\rad}_{L}\geq - \Lambda^{2}$) and $\epsilon \in (0,1)$ consider the subset $\Omega_{\phi_{h}}(\epsilon)\subset N\times L$  given by $$\Omega_{\phi_{h}}(\epsilon)=\left\{ (x,y)\in  N\times L\colon \phi_{h}(\rho_{N}(x)) \leq \log (\rho_{L}(y)+1)^{1-\epsilon}\right\}.  $$ Here $\rho_{L}(y)={\rm dist}_{L}(y_0, y)$, $y_0\in L$.

\begin{definition}   An isometric immersion  $\varphi \colon M \to N\times L$ of a Riemannian manifold $M$ into the product $N\times L$ is said to be $\phi_{h}$-bounded if there exists a compact $K\subset M$ and $\epsilon \in (0,1)$ such that $\varphi (M\setminus K)\subset \Omega_{h}(\epsilon)$.
\end{definition}
\begin{remark}The class of $\phi$-bounded submanifolds contains the class of  cylindrically bounded submanifolds.
\end{remark}

 \subsection{ Curvature estimates for $\phi$-bounded submanifolds}\label{sec:curvEst}In this section,  we extend   the cylindrically bounded version of Jorge-Koutrofiotis's  Theorem, Thm.\,\ref{thmABM}-ii. and the mean curvature estimates of  Thm.\,\ref{thm-Alias-Bessa-Dajczer} to the class of $\phi_{h}$-bounded properly immersed submanifolds. These extensions are done in two ways.  First:  the  class we consider is larger than the class of cylindrically bounded  submanifolds. Second:  there are no  requirements on the growth    on the second fundamental form as  in Thm.\,\ref{thmABM}. We also should observe that although  $\phi$-bounded properly immersed submanifolds,  ($\varphi \colon M \to N\times L$) are stochastically complete, provided $L$ has an Omori-Yau pair, see Section \ref{sec:OM-Pairs}, we do not need that to prove the following result.
\begin{theorem}
\label{thmMain-intro} Let  $\varphi\colon M\rightarrow N^{n-\ell}\times L^{\ell}$ be  a $\phi_{h}$-bounded isometric immersion of  a complete Riemannian $m$-manifold $M$ with  $n-m\leq m-\ell -1$. If $\varphi$ is proper and  $- G\leq b\leq 0$
then
\begin{equation}
sup_{M}K_{M}\geq \vert b\vert +\inf_{N}K_{N}.
\label{eq-ABM-2}
\end{equation} With strict inequality $\sup_{M}K_{M}>\inf_{N}K_{N}$ if $b=0$.
\end{theorem}

  \begin{corollary}
\label{CorMain-2} Let $\varphi\colon M\rightarrow  N^{n-\ell}\times L^{\ell}$ be a properly immersed, cylindrically bounded submanifold, $\varphi (M)\subset B_{N}(r)\times L^{\ell}$, where $B_{N}(r)$ is a regular geodesic  ball of $N$.   Suppose that $n-m \leq m-\ell-1$.
Then the sectional curvature of $M$ satisfies the following inequality
\begin{equation}
\sup_{M}K_{M}\geq C_{b}^2(r) +\inf_{N}K_{N},
\label{eq-ABM-4}
\end{equation}where  $b=\sup_{ B_{N}(r)}K_{N}^{rad}$ and $C_{b}$ is defined in \eqref{eqCb}.
\end{corollary}

Our next main result extends the mean curvature estimates \eqref{eqABD-mean} to $\phi$-bounded submanifolds.
\begin{theorem}Let  $\varphi\colon M\rightarrow N^{n-\ell}\times L^{\ell}$ be  a $\phi_{h}$-bounded isometric immersion of  a complete Riemannian $m$-manifold $M$ with  $m\geq \ell+1$. If $\varphi$ is proper then the mean curvature vector $H={\rm tr}\, \alpha$ of $\varphi$ satisfies\label{thmMain-2}
\begin{equation}\sup_{M}\vert H \vert \geq  (m-\ell)\cdot\inf_{r\in [0, \infty)} \frac{h'}{h}(r)\cdot \end{equation} If   $-G\leq b\leq 0$ then
\begin{equation}\sup_{M}\vert H \vert \geq  (m-\ell)\cdot\sqrt{\vert b\vert}. \end{equation}With strict inequality $\sup_{M}\vert H \vert >0$ if $b=0$.
\end{theorem}
\section{Proof of the main results}
\subsection{Basic results}

Let $M$ and $W$ be  Riemannian manifolds of dimension $m$ and $n$ respectively  and  let
 $\varphi\colon M\to W$ be an isometric immersion. For a given
 function
$g\in C^\infty(W)$  set $f=g\circ\varphi\in C^\infty(M)$.
Since
\[
\langle\grad_{\!_M}f,X\rangle=\langle\grad_{\!_W}g, X\rangle
\]
for every vector field $X\in TM$, we obtain
\[
\grad_{\!_W}g=\grad_{\!_M}f +(\grad_{\!_W}g)^{\perp}
\]
according to the decomposition $TW=TM\oplus T^\perp M$. An easy computation using the Gauss formula gives the
well-known relation (see e.g. \cite{jorge-koutrofiotis})
\begin{equation}
\label{eqBF2}
\HessM f(X,Y)= \Hess_{\!_W}g(X,Y) +\langle\grad_{\!_W}g,\alpha(X,Y)\rangle
\end{equation}
for all vector fields $X,Y\in TM$, where $\alpha$ stands for the second fundamental form of
$\varphi$. In particular, taking traces  with respect to an orthonormal frame
$\{ e_{1},\ldots, e_{m}\}$ in $TM$ yields
\begin{equation}\label{eqBF3}
\triangle_{_M}f =\sum_{i=1}^{m}\Hess_{\!_W}g(e_i,e_i)+ \langle\grad_{\!_W}g,H\rangle.
\end{equation}
where $H=\sum_{i=1}^{m}\alpha(e_i,e_i)$.

\vspace{2mm}
In the sequel, we will need  the following  well known results, see the classical Greene-Wu \cite{GW} for the Hessian Comparison Theorem and Pigola-Rigoli-Setti's \textquotedblleft must looking at\textquotedblright book \cite[Lemma 2.13]{prs-vanishing}, see also \cite{swanson}, \cite[Thm.1.9]{bianchini-mari-rigoli} for the Sturm Comparison Theorem.
\begin{theorem}[Hessian Comparison Thm.]\label{HCT} Let $W$ be a complete $n$-manifold and  $\rho_{_W}(x)={\rm dist}_{_W}(x_0,x)$, $x_0\in W$ fixed. Let $D_{x_{0}}=W\setminus (\{x_0\}\cup {\rm cut}(x_0))$ be the domain of normal geodesic coordinates at $x_0$. Let $G\in C^{0}([0, \infty))$  and let $h$ be the solution of  \eqref{eqg}. Let $[0, R)$ be the largest interval where $h>0$. Then
\begin{itemize}\item[i.]If the radial sectional curvatures along the geodesics issuing from $x_0$ satisfies $$K_{_W}^{rad}\geq -G(\rho_{_W}),\,\,{\rm in} \,\,B_{_W}(R) $$ then $$ \Hess_{\!_W}\rho \leq \frac{h'}{h}(\rho_{_W})\left[  \langle, \rangle - \di \rho \otimes \di \rho \right]\,\,{\rm on}\,\,D_{x_0}\cap B_{_W}(R)$$
\item[]
\item[ii.]If the radial sectional curvatures along the geodesics issuing from $x_0$ satisfy $$K_{_W}^{rad}\leq -G(\rho_{_W}),\,\,{\rm in} \,\,B_{_W}(R) $$ then $$ \Hess_{\!_W}\rho_{_W} \geq \frac{h'}{h}(\rho)\left[ \langle, \rangle - \di \rho \otimes \di \rho \right]\,\,{\rm on}\,\,D_{x_0}\cap B_{_W}(R)$$
\end{itemize}
\end{theorem}
\begin{lemma}[Sturm Comparison Thm.]\label{sturm}Let $G_1, G_2\in L^{\infty}_{loc}(\mathbb{R})$, $G_1\leq G_2$ and $h_1$ and $h_2$ solutions of the following problems:
\begin{equation}
\begin{array}{ll}a.)\left\{ \begin{array}{rll}   h_1''(t)-G_1 (t)h_1 (t) &\leq&0 \\
   h_1(0)= 0,\,\,\, \,\,
   h'_1(0)&>&0\end{array}\right. &\,\,\, b.)\left\{ \begin{array}{rll} h_2''(t)-G_2 (t)h_2 (t) &\geq &0 \\
   h_2(0)=0,\,\,\,\,\,
   h'_2(0)&>&h'_1(0),\end{array}\right.
\end{array}
\end{equation}and let $I_1=(0, S_1)$ and $I_{2}=(0,S_2)$ be the largest connected intervals where $h_1>0$ and $h_2>0$ respectively. Then
\begin{itemize}\item[1.] $S_1\leq S_2.$ \hspace{2mm}And on  $ I_1$, $\displaystyle\frac{h'_1}{h_1}\leq \frac{h'_2}{h_2}$ and  $h_1\leq h_2$.
\item[]
\item[2.] If $h_1(t_o)=h_2(t_o)$, $t_o\in I_1$ then $h_1\equiv h_2$ on $(0,t_o)$.
\end{itemize}
\end{lemma}
For a more detailed Sturm Comparison Theorem one should consult the beautiful book \cite[Chapter 2.]{prs-vanishing}.
If $-G=b\in \mathbb{R}$ then the solution of $h_{b}''(t)-G\cdot h_{b} (t)=0$ with $h_{b}(0)=0$ and $h_{b}'(0)=1$ is given by
 $$h_{b}(t)=\left\{\begin{array}{clrll}\displaystyle \frac{1}{\sqrt{-b}}\cdot \sinh( \sqrt{-b}\,t) &{\rm if}& b&<&0\\
 &&&&\\

  t &{\rm if}& b&=&0 \\
 &&&&\\ \displaystyle\frac{1}{\sqrt{b}}\cdot\sin(\sqrt{b}\,t) &{\rm if}& b&>&0.\end{array}\right.$$
In particular,  if  the radial sectional curvatures along the geodesics issuing from $x_0$ satisfy $K_{_W}^{rad}(x)\leq- G(\rho_{_W}(x))\leq  b$,  $x\in B_{_W}(R)=\{x,{\rm dist}_{_W}(x_0,x)= \rho_{_W}(x)<R\}$,  then the solution $h$  of \eqref{eqg}, satisfies
 $\displaystyle (h'/h)(t)\geq (h_{b}'/h_{b})(t)=C_{b}(t),$ $t\in (0, R)$,   $R< \pi/2\sqrt{b}$, if $b>0$.
 Therefore, $\Hess_{_W}\rho_{_W}\geq C_{b}(\rho_{_W})\left[\langle, \rangle -\di\rho_{_W}\oplus\di\rho_{_W}\right]$.
 Likewise, if  $K_{_W}^{rad}(x)\geq -G(\rho_{_W}(x))\geq b$, $x\in B_{_W}(R)$
  then $\displaystyle (h'/h)(t)\leq C_{b}(t),$ $t\in (0, R)$ and
   $\Hess_{_W}\rho_{_W}\leq C_{b}(\rho_{_W})\left[\langle, \rangle -\di\rho_{_W}\oplus\di\rho_{_W}\right]$.

\subsection{Proof of Theorem \ref{thmMain-intro}.} Assume without loss that there exists a $x_0\in M$ such that $\varphi(x_0)
  =(z_0,y_0)\in N\times L$, $z_0$, $y_0$ the distinguished points  of $N$ and $L$. For each  $x\in M$, let  $\varphi (x)=(z(x), y(x))$.
 Define $ g \colon N\times L \rightarrow \mathbb{R} $  by $ g(z,y) = \phi_{h}(\rho_{_N}(z))+1 $, recalling that $\phi_{h}(t)=\int_{0}^{t}h(s)ds$, and define   $ f=g\circ \varphi \colon M \rightarrow \mathbb{R} $  by $ f(x) = g(\varphi(x))=\phi_{h}(\rho_{_N}(z(x)))+1$. For each $k\in  \mathbb{N}$, set
 $g_{k}(x) = f(x) - \frac{1}{k}\cdot \log(\rho_{_L}(y(x))+1)$. Observe that $g_{k}(x_0)=1$ for all $k$, since $\rho_{_N}(z_0)= \rho_{_L}(y_0)=0$. First,  let us  prove the item i.

If $x \to \infty$ in $M$   then $\varphi (x)\to \infty$ in $N\times L$ since $ \varphi $ is proper. On the other hand,  $\varphi (M\setminus K)\subset \Omega_{h}(\epsilon)$ for some compact $K\subset M$ and $\epsilon \in (0,1)$. This implies that $y(x)\to \infty$ in $L$ and $$\frac{g_{k}(x)}{\log(\rho_{_L}(y(x))+1)}=\frac{f(x)}{\log(\rho_{_L}(y(x))+1)}-\frac{1}{k} <\frac{1}{\log(\rho_{_L}(y(x))+1)^{\epsilon}}-\frac{1}{k}<0$$for  $\rho_{_M}(x)\gg 1.$ This implies that $g_{k}(x)<0$  for $\rho_{_M}(x)\gg 1$. Therefore each $g_{k}$ reach a maximum at a point $x_{k}\in M$.  This  yields a sequence $ \{x_{k}\}\subset M$ so that $\Hess_{\!_M}g_{k}(x_k)(X,X)\leq 0$ for all $X\in T_{x_k}M$, this is, $\forall X\in T_{x_{k}}M $ \begin{equation}\HessM f(x_{k})(X,X)\leq \frac{1}{k}\cdot \HessM \log(\rho_{_L}( y(x_{k}))+1)(X,X).\label{eqHessiano-1}\end{equation}

 \noindent Observe  that $\log(\rho_{_L}( y(x_{k}))+1)= \log(\rho_{_L}\circ \pi_{L} +1) \circ \varphi (x_k)$, $\pi_{_L}\colon N\times L \to L$ the projection on the second factor,  thus the right hand side of \eqref{eqHessiano-1}, using the formula \eqref{eqBF2}, is given by
\begin{eqnarray}\HessM \log(\rho_{_L}( y(x_{k}))+1)(X,X)&=&\Hess_{\!_{N\times L}}\log(\rho_{_L}\circ \pi_{L} +1)(\varphi (x_k))(X,X)\nonumber \\
&&\label{eqHessiano2}\\
&+& \langle \grad_{\!_{N\times L}} \log(\rho_{_L}\circ \pi_{L} +1), \alpha(X,X) \rangle\nonumber\end{eqnarray}Where $\alpha$ is the second fundamental form of  $\varphi$.
For simplicity, set $\psi(t)=\log (t+1)$, $z_k=z(x_k)$,
 $y_{k}=y(x_{k})$, $s_k=\rho_{_N}(z_k)$ and $t_{k}=\rho_{_L}(y_k)$. Decomposing  $X\in TM$ as $X=X^{N}+X^{L}\in TN\oplus TL$, we see that the first term of the  right hand side  of  \eqref{eqHessiano2} is
\begin{eqnarray}
\Hess_{\!_{N\times L}}\psi\circ \rho_{_L}\circ y(x_{k})(X,X) &=&  \psi''(t_k)\vert X^{L}\vert^{2}+\psi'(t_k)\Hess_{\!_L} \rho_{_L}(y_k)(X,X)\nonumber\\
&&\nonumber\label{eqHessiano3-1} \\
&\leq &   \psi''(t_k)\vert X^{L}\vert^{2}+ C_{-\Lambda^{2}}(t_k)\frac{\vert X^{N}\vert^{2}}{(t_k+1)}\\
&&\nonumber \\
&\leq & C_{-\Lambda^{2}}(t_k)\frac{\vert X^{N}\vert^{2}}{(t_k+1)},\nonumber
\end{eqnarray}since $\Hess_{\!_L} \rho_{_L}(y_k)(X,X)\leq C_{-\Lambda^{2}}(t_k)\vert X^{N}\vert^{2} $ (by Theorem \ref{HCT}) and  $\psi'' \leq 0$.
\vspace{2mm}

 The second term of the  right hand side of   \eqref{eqHessiano2} is
 \begin{eqnarray}\langle \grad_{\!_{N\times L}} \psi\circ \rho_{_L}\circ y(x_{k}), \alpha(X,X)\rangle &=& \psi'(t_k)\langle \grad_{\!_L} \rho_{_L}(y_k), \alpha (X, X)\rangle  \nonumber\\
 && \nonumber \\
 &\leq & \frac{1}{(t_k+1)}\Vert \alpha\Vert \cdot \vert X  \vert^{2}  \label{eqHessiano4-1} \end{eqnarray}

 From \eqref{eqHessiano3-1} and \eqref{eqHessiano4-1} we have  the following
 \begin{eqnarray}\Hess_{\!_M}\psi\circ \rho_{_L}\circ y(x_{k})(X,X)&\leq & \frac{C_{-\Lambda^2}(t_k)+ \Vert \alpha\Vert  }{(t_k+1)}\cdot \vert X\vert^{2}\label{eqHessiano5-1}
 \end{eqnarray}And from \eqref{eqHessiano-1} and \eqref{eqHessiano5-1} we have that \begin{eqnarray}\Hess_{\!_M}f(x_{k})(X,X)& \leq &
 \frac{1}{k}\frac{(C_{-\Lambda^2}(t_k)+ \Vert \alpha\Vert)}{(t_k+1)}\vert X\vert^{2}\label{eqHessianoI-1}
 \end{eqnarray}

We will compute the left hand side of \eqref{eqHessiano-1}. Using the formula \eqref{eqBF2} again we have
\begin{eqnarray}\Hess_{\!_M}f(x_{k}) &=& \Hess_{\!_{N\times L}} g(\varphi (x_k)) + \langle \grad_{\!_{N\times L}}\, g, \alpha \rangle\label{eqHessiano6-1}
\end{eqnarray}Recalling that $f=g\circ\varphi$ and $g$ is given by  $ g(z,y) = \phi_{h}(\rho_{_N}(z)) $, where $\phi_{h}$ is defined in \eqref{eq-phi}
and $ \rho_{_N}(z) = dist_{N}(z_{0},z) $.  Let us consider an orthonormal basis \eqref{base} \begin{equation}\{\stackrel{\in TN}{\overbrace{ \grad \rho_{\!_N}, \partial/\partial \theta_{\!_1}, \ldots, \partial/\partial \theta_{\!_{n-\ell-1}}}}, \stackrel{\in TL}{\overbrace{ \partial/\partial \gamma_{\!_1}, \ldots, \partial/\partial \gamma_{\!_\ell} }}\}\label{base}\end{equation} for $T_{\varphi (x_k)}(N\times L)$. Thus if $X\in T_{x_k}M$, $\vert X\vert =1$, we can decompose $$X=a\cdot \grad \rho_{_N}+\sum_{j=1}^{n-\ell-1}b_{j}\cdot \partial/\partial \theta_{j} + \sum_{i=1}^{\ell}c_{i}\cdot \partial/\partial \gamma_i$$  with $a^2+\sum_{j=1}^{n-\ell-1} b_{j}^{2}+ \sum_{i=1}^{\ell} c_{i}^{2}=1$. Recalling that $s_k=\rho_{_N}(z(x_k))$, we can see that   the first term of the right hand side of \eqref{eqHessiano6-1}
\begin{eqnarray} \Hess_{\!_{N\times L}} g(\varphi (x))(X,X)\!&=&\!\phi_{h}''(s_k)\!\cdot  \! a^2+ \phi_{h}'(s_k)\!\sum_{j=1}^{n-\ell-1}\!b_j^{2}\cdot \!\Hess \rho_{_N}\!(z_k)(\frac{\partial}{\partial \theta_{j}},\frac{\partial}{\partial \theta_{j}})\nonumber\\
&&\nonumber \\
&\geq &\!\phi_{h}''(s_k)\cdot   a^2+ \phi_{h}'(s_k)\sum_{j=1}^{n-\ell-1}\!b_j^{2}\cdot \frac{h'}{h}(s_k)\nonumber\\
&&\nonumber \\
&=& \phi_{h}''(s_k)\cdot a^2+   (1-a^2 - \sum_{i=1}^{\ell}\!c_i^{2})\cdot \phi_{h}'(s_k)\cdot \frac{h'}{h}(s_k)\nonumber \\
&& \nonumber \\
&=& \left[\stackrel{\equiv 0}{(\overbrace{\phi_{h}''-\frac{h'}{h}\cdot \phi_{h}'})}a^2 + (1- \sum_{i=1}^{\ell}\!c_i^{2})\cdot \phi_{h}'\cdot \frac{h'}{h}\right](s_k) \nonumber\\
&&\nonumber \\
&=& (1- \sum_{i=1}^{\ell}\!c_i^{2})\cdot \phi_{h}'(s_k)\cdot \frac{h'}{h}(s_k)\nonumber
\end{eqnarray}

Thus \begin{equation}\label{eqHessiano7}\Hess_{\!_{N\times L}} g(\varphi (x))(X,X)\geq (1- \sum_{i=1}^{\ell}\!c_i^{2})\cdot \phi_{h}'(s_k)\cdot \frac{h'}{h}(s_k).\end{equation}

The second term of the right hand side of \eqref{eqHessiano6-1} is the following
\begin{eqnarray} \langle \grad_{\!_{N\times L}} g, \alpha (X,X)\rangle & = & \phi_{h}'(s_k)\langle \grad_{\!_N} \rho_{_N}(z_k), \alpha (X,X)\rangle\nonumber\\
&& \label{eqHessiano8}\\
&\geq & -\phi_{h}'(s_k)\vert \alpha(X,X) \vert\nonumber
\end{eqnarray}

From \eqref{eqHessiano6-1}, \eqref{eqHessiano7}, \eqref{eqHessiano8} we have that,

\begin{equation}\Hess_{\!_M}f(x_{k})(X,X) \geq \left[(1- \sum_{i=1}^{\ell}\!c_i^{2})\cdot  \frac{h'}{h}(s_k) -  \vert \alpha (X ,X) \vert\right]\phi_{b}'(s_k)\label{eqHessiano9-1}\end{equation}

Recall that $n+\ell\leq 2m-1$. This dimensional restriction implies that $m\geq \ell + 2$, since $n\geq m+1$. Therefore, for every $x\in M$ there exists a  sub-space $V_x\subset T_xM$ with ${\rm dim}(V_x)\geq (m-\ell)\geq  2$ such that $V\perp TL$, this is equivalent to $c_i=0$. If we take any $X\in V_{x_k}\subset T_{x_k}M$, $\vert X\vert=1$ we have by \eqref{eqHessiano9-1} that

\begin{equation*} \frac{(C_{-\Lambda^2}(t_k)+ \vert \alpha(X,X)\vert)}{k(t_k+1)}\geq \Hess_{\!_M}f(x_{k})(X,X) \geq \left[ \frac{h'}{h}(s_k) -  \vert \alpha (X ,X) \vert\right]\phi_{h}'(s_k)\label{eqHessiano10}\end{equation*}
Thus, reminding that $\phi'_h=h$,  \begin{eqnarray} \label{eqHessiano11-p} \left[\frac{1}{k(t_k+1)} + h(s_k)\right]\vert \alpha (X ,X) \vert\geq h'(s_k)-\frac{C_{-\Lambda^{2}}(t_k)}{k(t_k+1)}
\end{eqnarray}

Since $-G\leq b \leq 0$,  we have by Lemma \ref{sturm} (Sturm's argument) that the solution $h$ of \eqref{eqg} satisfies $ \displaystyle (h'/h)(t)\geq C_{b}(t)>\sqrt{\vert b\vert}$ and that $h(t)\to +\infty$ as $t\to +\infty$, where $C_b$ is defined in \eqref{eqCb}.
Let us assume  that $x_{k}\to \infty$ in $M$, (the case $\rho_{M}(x_k)\leq C^{2} < \infty$ will be considered later), then $s_{k}\to \infty$ as well as $t_{k}\to \infty$. Thus from \eqref{eqHessiano11-p},  for sufficiently  large $k$,  we have  at $\varphi (x_k)$ that
\begin{eqnarray} \label{eqHessiano11-b} \left[\frac{1}{k(t_k+1)h(s_k)} + 1\right]\vert \alpha (X ,X) \vert &\geq & \frac{h'(s_k)}{h(s_k)}-\frac{C_{-\Lambda^{2}}(t_k)}{k(t_k+1)h(s_k)}\nonumber \\
&\geq & C_{b}(s_k)-\frac{C_{-\Lambda^{2}}(t_k)}{k(t_k+1)h(s_k)}\nonumber\\
&>& 0
\end{eqnarray} Thus, at $x_k$ and $X\in T_{x_k}M$ with $\vert X\vert=1$ we have
\begin{equation}\label{eq27}\vert\alpha (X,X)\vert \geq\left[ C_{b}(s_k)-\frac{C_{-\Lambda^{2}}(t_k)}{k(t_k+1)h(s_k)}\right]\left[\frac{1}{k(t_k+1)h(s_k)} + 1\right]^{-1}>0.\end{equation}  We will need the following lemma known as Otsuki's Lemma \cite[p.28]{kobayashi-nomizu}. \begin{lemma}[Otsuki] Let $ \beta \colon \mathbb{R}^{q}\times\mathbb{R}^q \rightarrow \mathbb{R}^d $, $ d \leq q-1 $, be a symmetric bilinear form satisfying $ \beta(X,X)\neq 0 $ for $ X\neq 0 $. Then there exists linearly independent vectors $ X,Y $ such that $\beta(X,X)=\beta(Y,Y) $ and $ \beta(X,Y)=0 $.
\end{lemma}\noindent The {\em horizontal} subspace $V_{x_k}$ has dimension ${\rm dim}(V_{x_k})\geq m-\ell\geq 2$. Thus, by the inequality \eqref{eq27} and  $n-m\leq m-\ell-1\leq {\rm dim}(V_{x_k})-1 $,  we may apply Otsuki's Lemma to   $\alpha (x_k)\colon  V_{x_k}\times V_{x_k}\to T_{x_k}M^{\perp}\simeq \mathbb{R}^{n-m}$ to  obtain $X, Y\in V_{x_k}$, $\vert X\vert \geq \vert Y\vert \geq 1$ such that $\alpha (x_k)(X,X)= \alpha (x_k)(Y,Y)$ and $\alpha (x_k)(X,Y)=0$.

\vspace{2mm}

By the Gauss equation we have  that
\begin{eqnarray}\label{curvaturek}
K_{M}(x_k)(X,Y) - K_{N}(\varphi(x_k))(X,Y) &=& \frac{\langle \alpha (x_k)(X,X),\alpha (x_k)(Y,Y)\rangle  }{\vert X\vert^{2}\vert Y\vert^{2} - \langle X,Y\rangle^{2}} \nonumber \\
&=& \frac{\vert \alpha (x_k)(X,X)\vert^{2}}{\vert X\vert^{2}\vert Y\vert^{2}} \nonumber \\
&\geq & \left(\frac{\vert \alpha (x_k)(X,X)\vert }{\vert X\vert^{2}}\right)^{2}\nonumber \\
&=& \left\vert \alpha (x_k)\left(\frac{X}{\vert X\vert },\frac{X}{\vert X\vert }\right)\right\vert^{2}\nonumber
\end{eqnarray}This implies by \eqref{eq27} that $$\sup K_{M}-\inf K_{N}>\left( \left[\frac{h'(s_k)}{h(s_k)}-\frac{C_{-\Lambda^{2}}(t_k)}{k(t_k+1)h(s_k)}\right]\left[\frac{1}{k(t_k+1)h(s_k)} + 1\right]^{-1}\right)^{2}>0.$$ Therefore,  $ \sup K_{M}-\inf K_{N}>0$ regardless $b=0$ or $b<0$.  If $b<0$ we let $k \to +\infty$ and then we have
\begin{eqnarray}\sup K_{M}-\inf K_{N}&\geq& \lim_{s_{k}\to \infty} \left[\frac{h'}{h} (s_k)\right]^2= \vert b\vert
\end{eqnarray}

The  case where the sequence $ \{x_k\} \subset M $ remains in a compact set, we proceed as follows.  Passing to a subsequence we have that $ x_k \rightarrow x_{\infty} \in M $. Thus $t_{k}\to t_{\infty}<\infty$ and $s_{k}\to s_{\infty}<\infty$. By \eqref{eqHessianoI-1}
\begin{eqnarray}\Hess_{M}f(x_{\infty})(X,X)& \leq &
 \lim_{k\to\infty} \frac{(C_{-\Lambda^{2}}(t_{\infty})+  \vert \alpha (x_{\infty})(X,X)\vert )}{k(t_{\infty}+1)}=0,\label{eqHessianoI-I}
 \end{eqnarray}
 for all $ X \in T_{x_0}M $. Using the expression on the right hand side of \eqref{eqHessiano9-1} we obtain for every $ X \in V_{x_\infty} $
\begin{eqnarray*}
0 \geq \Hess f(x_\infty)(X,X) \geq \left[(1- \sum_{i=1}^{\ell}\!c_i^{2})\cdot  \frac{h'}{h}(s_\infty) -  \vert \alpha (X ,X) \vert\right]\phi_{b}'(s_\infty) .
\end{eqnarray*}There exists a  sub-space $V_x\subset T_xM$ with ${\rm dim}(V_x)\geq (m-\ell)\geq  2$ such that $V\perp T\mathbb{R}^{\ell}$, this is equivalent to $c_i=0$. If we take any $X\in V_{x_\infty}\subset T_{x_\infty}M$, $\vert X\vert=1$ we have hence
\begin{eqnarray*}
\vert \alpha_{x_\infty}(X,X)\vert  \geq \frac{h'}{h}(s_{\infty})\vert X\vert^2 .
\end{eqnarray*}

Again, using Otsuki's Lemma and Gauss equation, we conclude that
\begin{eqnarray}
\sup_{M}K_{M} - \inf_{B_{N}(r)}K_{N} \geq \frac{h'}{h}(s_{\infty}) > \vert b\vert.
\end{eqnarray}

\subsection{Proof of Theorem \ref{thmMain-2}. }We will follow the proof of Theorem \ref{thmMain-intro} closely. Recall  that $g_{k}$ reaches a maximum at $ x_{k}\in M$, $k=1,2,\ldots$, thus  so that $\triangle_{\!_M}g_{k}(x_k)\leq 0$. Thus  \begin{equation}\triangle_{_M}f(x_{k})\leq \frac{1}{k}\cdot \triangle_{_M}( \log(\rho_{_L}\circ \pi_{L} +1) \circ \varphi (x_k)).\label{eqLaplacianoP}\end{equation} Using the formula \eqref{eqBF3}
\begin{eqnarray} \triangle_{_M}( \log(\rho_{_L}\circ \pi_{L} +1) \circ \varphi (x_k))&=&\sum_{i=1}^{m}\Hess_{N\times L}\log(\rho_{\mathbb{R}^{\ell}}\circ \pi_{L} +1)(\varphi (x_k))(X_{i},X_{i})\nonumber \\
&&\label{eqHessiano2A}\\
&&+\, \langle \grad_{_{N\times L}} \log(\rho_{_L}\circ \pi_{L} +1), H \rangle\nonumber\end{eqnarray}where $H=\sum_{i=1}^{m}\alpha(X_i,X_i)$ is the mean curvature vector while $\alpha$ is the second fundamental form of the immersion $\varphi$ and $\{X_i\}$ is an orthonormal basis of $T_{x_k}M$.
\vspace{2mm}

As before,  decomposing  $X\in TM$ as $X=X^{N}+X^{L}\in TN\oplus TL$  and setting $\psi(t)=\log (t+1)$,
 $y_{k}=y(x_{k})$ and $t_{k}=\rho_{_L}(y_k)$ we have that the right hand side of  \eqref{eqHessiano2A}
\begin{eqnarray}
\sum_{i=1}^{m}\Hess_{N\times L}\psi\circ \rho_{_L}\circ y(x_{k})(X_i,X_i) &=&  \psi''(t_k)\sum_{i=1}^{m}\vert X_{i}^{L}\vert^{2}\nonumber \\ && +\,\psi'(t_k)\sum_{i=1}^{m}\Hess_{L} \rho_{_L}(y_k)(X_i,X_i)\nonumber\\
&&\label{eqHessiano3} \\
&\leq & \frac{C_{-\Lambda^{2}}(t_k)}{(t_k+1)}\sum_{i=1}^{m}\vert X_{i}^{N}\vert^{2},\nonumber\\
&& \nonumber \\
& \leq & \frac{m\cdot  C_{-\Lambda^{2}}(t_k)}{(t_k+1)}\nonumber
\end{eqnarray}since $\psi'' \leq 0$ and
 \begin{eqnarray}\langle \grad_{N\times L} \psi\circ \rho_{_L}\circ y(x_{k}), H\rangle &=& \psi'(t_k)\langle \grad \rho_{_L}(y_k), H\rangle  \nonumber\\
 && \nonumber \\
 &\leq & \frac{1}{(t_k+1)}\vert H\vert  \label{eqHessiano4} \end{eqnarray}

 From \eqref{eqHessiano2A}, \eqref{eqHessiano3} and \eqref{eqHessiano4} we have
 \begin{eqnarray}\triangle_{M} \log(\rho_{_L}( y(x_{k}))+1)&\leq & \frac{m\cdot  C_{-\Lambda^{2}}(t_k)+ \vert H\vert  }{(t_k+1)}\label{eqHessiano5}
 \end{eqnarray}And from \eqref{eqLaplacianoP} and \eqref{eqHessiano5} we have that \begin{eqnarray}\triangle_{M}f(x_{k})& \leq &
 \frac{m\cdot  C_{-\Lambda^{2}}(t_k)+ \vert H\vert  }{k(t_k+1)}\label{eqHessianoI}
 \end{eqnarray}

We will compute the left hand side of \eqref{eqLaplacianoP}. Recall that $f=g\circ\varphi$ and $g$ is given by  $ g(z,y) = \phi_{h}(\rho_{_N}(z)) $, where $\phi_{}$ is defined in \eqref{eq-phi}. Using the formula \eqref{eqBF3} again we have
\begin{eqnarray}\triangle_{M}f(x_{k})&=& \sum_{i=1}^{m}\Hess_{N\times L} g(\varphi (x_k))(X_i,X_i) + \langle \grad\, g, H\rangle\label{eqHessiano28}
\end{eqnarray}  Consider the  orthonormal basis  \eqref{base} for $T_{\varphi (x_k)}(N\times L)$. Thus if $X_i\in T_{x_k}M$, $\vert X_i\vert =1$, we can decompose $$X_i=a_i\cdot \grad \rho_{_N}+\sum_{j=1}^{n-\ell-1}b_{ij}\cdot \partial/\partial \theta_{j} + \sum_{l=1}^{\ell}c_{il}\cdot \partial/\partial \gamma_l$$  with $a_i^2+\sum_{j=1}^{n-\ell-1} b_{ij}^{2}+ \sum_{l=1}^{\ell} c_{il}^{2}=1$. Set $z_k=z(x_k)$ and $s_k=\rho_{_N}(z_k)$. We have  as in \eqref{eqHessiano7}
\begin{eqnarray} \Hess_{N\times L} g(\varphi (x))(X_i,X_i)&\geq &(1- \sum_{l=1}^{\ell}\!c_{il}^{2})\cdot \phi_{h}'(s_k)\cdot \frac{h'}{h}(s_k)\label{eqHessiano29}
\end{eqnarray}

\noindent The second term of the right hand side of \eqref{eqHessiano28} is the following, if $\vert X\vert =1$,
\vspace{2mm}
\begin{eqnarray} \langle \grad g, H\rangle & = & \phi_{h}'(s_k)\langle \grad \rho_{_N}(z_k), H\rangle\nonumber\\
&& \label{eqHessiano30}\\[-1mm]
&\geq & -\phi_{h}'(s_k)\vert H \vert\nonumber
\end{eqnarray}
\vspace{1mm}

Therefore from \eqref{eqHessiano28}, \eqref{eqHessiano29}, \eqref{eqHessiano30} we have that,

\begin{equation}\triangle_{M}f(x_{k}) \geq \left[(m- \sum_{i=1}^{m}\sum_{l=1}^{\ell}\!c_{il}^{2})\cdot  \frac{h'}{h}(s_k) -  \vert H \vert\right]\phi_{b}'(s_k)\label{eqHessiano31}\end{equation}

From \eqref{eqHessianoI} and \eqref{eqHessiano31} we have

\begin{equation}\frac{m\cdot  C_{-\Lambda^{2}}(t_k)+ \vert H\vert  }{k(t_k+1)}\geq \triangle_{M}f(x_{k}) \geq \left[ (m-\ell)\cdot\frac{h'}{h}(s_k) -  \vert H \vert\right]\phi_{h}'(s_k)\label{eqHessiano10-B}\end{equation}
Therefore
 \begin{eqnarray} \label{eqHessiano33} \sup_{M}\vert H \vert\left[\frac{1}{h(s_k)\cdot k \cdot (t_k+1)} + 1\right] & \geq &  (m-\ell)\cdot\frac{h'}{h}(s_k) -\frac{m\cdot  C_{-\Lambda^{2}}(t_k)}{h(s_k)\cdot k \cdot (t_k+1)}\nonumber
\end{eqnarray}
Letting  $k\to \infty$ we have   $$\sup_{M}\vert H \vert \geq  (m-\ell)\cdot\lim_{k\to \infty} \frac{h'}{h}(s_k)\cdot$$

If in addition, we have that $-G\leq b \leq 0$ then  $ \displaystyle (h'/h)(s)\geq C_{b}(s)$. The case that $b=0$ we have $\displaystyle (h'/h)(s_k) \geq 1/s_k$ and $h(s_k)\geq s_k$. Since the immersion is $\phi$-bounded we have  $s_{k}^{2} \leq 2\log (t_{k}+1)^{(1-\epsilon)} $. Thus for sufficient  large $k$
$$ \sup_{M}\vert H \vert\left[\frac{1}{s_k\cdot k \cdot (t_k+1)} + 1\right]\geq \frac{ m-\ell}{s_k} -\frac{m\cdot  C_{-\Lambda^{2}}(t_k)}{s_k\cdot k \cdot (t_k+1)}>0.$$
This shows that $\sup_{M}\vert H\vert >0$.
\vspace{2mm}

In the case $b<0$, we have  $ (h'/h)(s_k)\geq C_{b}(s_k)\geq \sqrt{\vert b\vert}$ and $$\sup_{M}\vert H \vert \geq  (m-\ell)\cdot\lim_{k\to \infty} \frac{h'}{h}(s_k)\geq \sqrt{\vert b\vert}.$$
\begin{remark} The statements of Theorems \ref{thmMain-intro} and \ref{thmMain-2} are also true in a slightly more general situation. This is, if, instead a proper $\phi$-bounded immersion, one asks a proper immersion $\varphi\colon M\to N\times L$ with the property $$\lim_{x\to \infty_{in\,M}} \frac{\phi_{h}(\rho_{_N}(z(x)))}{\log (\rho_{_L}(y(x))+1)}=0 ,$$ where $\varphi(x)=(z(x),y(x))\in N\times L$.
\end{remark}

\section{ Omori-Yau pairs}\label{sec:OM-Pairs}

 Omori, in \cite{omori}, discovered an important global  maximum principle  for complete Riemannian manifolds with sectional curvature bounded below.  Omori's maximum principle was refined and extended by  Cheng and Yau, \cite{cheng-yau}, \cite{yau},  \cite{yau2},  to Riemannian manifolds with Ricci curvature bounded below and applied to find elegant solutions to various analytic-geometric problems on Riemannian manifolds.
There were others generalizations of the Omori-Yau maximum principle under more relaxed curvature requirements in  \cite{chen-xin},  \cite{dias} and an extension to an all general setting by  S. Pigola, M. Rigoli and  A. Setti in their beautiful book \cite{prs-memoirs}. There, they introduced the following terminology.
\begin{definition}[Pigola-Rigoli-Setti] The Omori-Yau maximum principle holds on a  Riemannian manifolds $W$ if for any  $u\in C^{2}(W)$ with $u^{\ast}\colon=\sup_{\!_W}u<\infty$, there exists a sequence of points $x_{k}\in W$, depending on $u$ and on $W$,  such that
\begin{eqnarray}\label{eq-omori-2}
\lim_{k\to \infty}u(x_{k}) =u^{\ast},& \vert \grad u\vert  (x_{k}) < \displaystyle \frac{1}{k},&  \triangle u(x_{k}) <  \displaystyle \frac{1}{k}.
\end{eqnarray} Likewise, the Omori-Yau maximum principle \textit{for the Hessian} holds on $W$ if  \begin{eqnarray}\label{eq-omori}
\lim_{k\to \infty}u(x_{k}) =u^{\ast},& \vert \grad u\vert  (x_{k}) < \displaystyle \frac{1}{k},&  \Hess_{\!_W} u(x_{k})(X,X) <  \displaystyle \frac{1}{k}\cdot \vert  X \vert^{2},
\end{eqnarray}for every $X\in T_{x_k}W$.
 \end{definition}A natural and important  question is, what are the Riemannian geometries  the Omori-Yau maximum principle holds on? It does hold on  complete Riemannian manifolds with sectional curvature bounded below holds \cite{omori}, it holds on  complete Riemannian manifolds with Ricci curvature bounded below \cite{cheng-yau}, \cite{yau},  \cite{yau2}. Follows from the work of Pigola-Rigoli-Setti \cite{prs-memoirs} that the Omori-Yau maximum principle  holds on complete Riemannian manifolds $W$ with  Ricci curvature with strong quadratic decay, $$ Ric_{_W} \geq - c^{2}\cdot\rho_{\!_W}^{2}\cdot\Pi_{i=1}^{k}(\log^{(i)}(\rho_{\!W}+1),\,\,\,\rho_{\!_W}\gg 1. $$ The notion of
 Omori-Yau pair was formalized in \cite{alias-bessa-montenegro-piccione}, after the work of Pigola-Rigoli-Setti. The Omori-Yau pair is,  here, described for the Laplacian and for the Hessian however, it certainly can be extended to other operators or bilinear forms.

\begin{definition}\label{def:OY-pair}Let $W$ be a  Riemannian manifold.
A pair  $ (\mathcal{G},\gamma)$ of smooth functions $ \mathcal{G} \colon [0,+\infty) \rightarrow (0,+\infty) $,
$ \gamma \colon W \rightarrow [0,+\infty) $, $\mathcal{G}\in C^{1}([0, \infty)),\,\gamma \in C^{2}([0, \infty))$,  forms an Omori-Yau pair for the Laplacian in W, if they satisfy the following conditions:\begin{itemize}
\item[h.1)] $\gamma (x) \rightarrow  +\infty $ as $x\rightarrow \infty$ ${\rm in\;W} $.
\item[]
\item[h.2)] $\mathcal{G}(0)>0$, $ \mathcal{G}'(t) \geq 0$ and  $\displaystyle{\int_{0}^{+\infty}\frac{ds}{\sqrt{\mathcal{G}(s)}} = +\infty}.$
\item[]
\item[h.3)]$\exists A>0$ constant such that $ \vert \grad_{_W} \gamma\vert \leq A  \sqrt{\mathcal{G}(\gamma)}\displaystyle\left(\int_{0}^{\gamma}\frac{ds}{\sqrt{\mathcal{G}(s)}}+1\right) $\, off a compact set.
\item[]
 \item[h.4)] $\exists B>0$  constant such that  $\triangle_{_W} \gamma  \leq B  \sqrt{\mathcal{G}(\gamma)}\displaystyle\left(\int_{0}^{\gamma}\frac{ds}{\sqrt{\mathcal{G}(s)}}+1\right)$\,\, off a compact set.
\end{itemize}
    The pair $ (\mathcal{G},\gamma)$ forms an Omori-Yau pair for the Hessian if instead h.4) one has
    \begin{itemize}\item[h.5)] $\exists C>0$  constant such that  $\Hess \gamma \leq C  \sqrt{\mathcal{G}(\gamma)}\displaystyle\left(\int_{0}^{\gamma}\frac{ds}{\sqrt{\mathcal{G}(s)}}+1\right) $  off a compact set,  in the sense of quadratic forms.
    \end{itemize}
\end{definition}

The result \cite[Thm.1.9]{prs-memoirs}  captured the essence of the Omori-Yau maximum principle and it
can be stated as follows.
\begin{theorem}If a Riemannian manifold $M$ has  an Omori-Yau pair $(\mathcal{G}, \gamma)$ then the Omori-Yau maximum principle on it.
\end{theorem}\noindent  The main step in the proof of  Alias-Bessa-Montenegro's Theorem (Thm.\ref{thmABM}) and Alias-Bessa-Dajczer's Theorem (Thm.\ref{thm-Alias-Bessa-Dajczer}) is to show   that a cylindrically bounded submanifold, properly immersed  into $ N\times L$,  with {\em  controlled} second fundamental form or {\em  controlled} mean curvature vector, has  an Omori-Yau pair, provided $L$ has an Omori-Yau pair. Thus,  the Omori-Yau maximum principle holds on those submanifolds and their proof follows the steps of Jorge-Koutrofiotis's Theorem.  On the other hand,  the idea behind  the proof of Theorems \ref{thmMain-intro} \& \ref{thmMain-2} is that: the factor $L$ has bounded sectional curvature it has a natural  Omori-Yau pair ($\mathcal{G}, \gamma$). This Omori-Yau pair together with  the geometry of the factor $N$  allows us to consider  an unbounded region $\Omega_{\phi}$  such that if $\varphi \colon M \to \Omega_{\phi} \subset N\times L$ is an isometric immersion then there exists a  function $f\in C^{2}(M)$, not necessarily bounded, and a sequence $x_k\in M$ satisfying $\triangle f (x_k)\leq 1/k$.  We show   that a properly immersed $\phi$-bounded submanifold has an Omori-Yau pair for the Laplacian, provided  the fiber $L$ has an Omori-Yau pair for the Hessian. We show in Theorem \ref{thm-OYBFL} that an Omori-Yau pair for the Hessian guarantee the Omori-Yau sequence for certain unbounded functions, as  this unbounded function $f$ we are working. This leads to  stronger  forms of Theorem \ref{thmMain-intro}. \& Theorem \ref{thmMain-2}.

Let $M$, $N$, $L$ be complete Riemannian manifolds of dimension $m$, $n-\ell$ and $\ell$,   with distinguished points $x_0$, $z_0$ and $y_0$  respectively. Suppose that $K_{N}^{\rad}\leq - G(\rho_{_N})$, $G$ satisfying \eqref{buonacurva2}. Let $h$ solution of \eqref{eqg} and $\phi_{h}$ as in \eqref{eq-phi}. Suppose in addition that $L$ has an Omori-Yau pair for the Hessian $(\gamma, \mathcal{G})$. Let $\Omega_{h, \gamma, \mathcal{G}}(\epsilon)\subset N\times L$ be the region defined by $$\Omega_{h, \gamma, \mathcal{G}}(\epsilon) =\{(z, y)\in N\times L\colon \phi_{h}\circ\rho_{_N}(z(x))\leq \left[\psi\circ\gamma (y(x))\right]^{1-\epsilon}\}, $$ where $ \displaystyle{\psi(t) = \log\left(\int_{0}^{t}\frac{ds}{\sqrt{\mathcal{G}(s)}}+1\right)}$. In this setting we have the following result.
\begin{theorem}\label{ThmBFP}
Let $ \varphi \colon M \to  N\times L$ be a properly immersed submanifold  such that $\varphi (M\setminus K )\subset \Omega_{h, \gamma, \mathcal{G}}(\epsilon) $ for some compact $K\subset M$ and positive $\epsilon \in (0,1)$. \begin{itemize}\item[1.]If
 $K_{N}^{\rad}\leq -G\leq b\leq 0$ and the codimension satisfies $n-m\leq m-\ell -1$ then \begin{equation}\sup_{M}K_{M}\geq \vert b\vert + \inf_{N}K_{N}.\end{equation} With strict inequality $\sup_{M}K_{M}>\inf_{N}K_{N}$ if $b=0$.
\item[2.]If $m\geq \ell +1$ then \begin{equation}\sup_{M}\vert H \vert \geq  (m-\ell)\cdot\inf_{r\in [0, \infty)} \frac{h'}{h}(r)\cdot \end{equation} If   $-G\leq b\leq 0$ then
\begin{equation}\sup_{M}\vert H \vert \geq  (m-\ell)\cdot\sqrt{\vert b\vert}. \end{equation}With strict inequality $\sup_{M}\vert H \vert >0$ if $b=0$.\end{itemize}

%
%
\end{theorem}
 Assume without loss of generality that there exists  $x_0\in M$ such that $\varphi(x_0)
  =(z_0,y_0)\in N\times L$. As before, $\varphi (x)=(z(x), y(x))$ and  $g,p\colon N\times L \to \mathbb{R}$ given by $g(z, y)=\phi_{h}(\rho_{N}(z)) + \psi(\gamma (y))$, $p(z,y)=\psi(\gamma(y))$.

  For each $k\in \mathbb{N}$, let $g_{k}\colon M \to \mathbb{R}$ given by $g_{k}(x)=g\circ \varphi (x) - p\circ \varphi(x)/k$. Observe that  $g_k(x_0)=1$ and for $\rho_{_M}(x)\gg 1$, we have that $g_{k}(x)<0$. This implies that $g_k$ has a maximum at a point $x_k$, yielding in this way a sequence $\{x_k\}\subset M$ such that $\HessM g_k(x_k)\leq 0$ in the sense of quadratic forms. Proceeding as in the proof of Theorem \ref{thmMain-intro} we have that for $X\in T_{x_k}M$,
  \begin{equation}\HessM g\circ \varphi (x_k)(X,X)\leq \frac{1}{k} \HessM p \circ \varphi(x_k)(X,X). \label{eqHessGP-1}\end{equation} We have to compute both terms of this inequality. Considering once more the orthonormal basis \eqref{base} for $T_{\varphi (x_k)}(N\times L)$ we can decompose,   $X\in T_{x_k}M$, $\vert X\vert =1$, (after identifying $X$ with $d\varphi X$), as $$X=a\cdot \grad \rho_{_N}+\sum_{j=1}^{n-\ell-1}b_{j}\cdot \partial/\partial \theta_{j} + \sum_{i=1}^{\ell}c_{i}\cdot \partial/\partial \gamma_i$$  with $a^2+\sum_{j=1}^{n-\ell-1} b_{j}^{2}+ \sum_{i=1}^{\ell} c_{i}^{2}=1$. Setting $s_k=\rho_{N}(z(x_k))$, $t_k=\gamma(y(x_k))$, we have as in \eqref{eqHessiano9-1},
   \begin{eqnarray}\HessM g\circ \varphi (x_k)(X,X)&=&\Hess_{\!_{N\times L}} g (\varphi(x_k))(X,X)+ \langle \grad_{\!_{N\times L}}g, \alpha(X,X)\rangle\nonumber \\ \label{eqHessGP-2}
  &\geq &\left[(1- \sum_{i=1}^{\ell}\!c_i^{2})\cdot  \frac{h'}{h}(s_k) -  \vert \alpha (X ,X) \vert\right]\phi_{b}'(s_k)\\
  && \nonumber  \end{eqnarray}
  \begin{eqnarray}
  \HessM p \circ \varphi(x_k)(X,X)&=&\Hess_{\!_{N\times L}} p (\varphi(x_k))(X,X)+ \langle \grad_{\!_{N\times L}}p, \alpha(X,X)\rangle\nonumber \\
   && \nonumber \\
  &=&\psi''(t_k)\langle X, \grad_{_L}\! \gamma  \rangle^{2} \! +\!\psi'(t_k)\Hess_{\!_L} \gamma (X, X)\nonumber \\
   && \nonumber \\
  & &+\,\,\psi'(t_k)\langle \grad_{\!_{L}}\gamma , \alpha(X,X)\rangle \nonumber\\
   && \nonumber \\
  &\leq & \psi'(t_k)\left(\Hess_{\!_L} \gamma (X, X) + \vert  \grad_{\!_{L}}\gamma\vert \cdot \vert \alpha(X,X)\vert \right) \label{eqHessGP-3} \\
   && \nonumber \\
  &\leq & \displaystyle\frac{\left[\sqrt{\mathcal{G}(\gamma (t_k))}\displaystyle\left(\int_{0}^{t_{k}} \frac{ds}{\sqrt{\mathcal{G}(\gamma (s))}}+1\right)\right]\left(C+ A \cdot \vert \alpha(X,X)\vert\right)}{\sqrt{\mathcal{G}(\gamma (t_k))}\displaystyle\left(\int_{0}^{t_{k}} \frac{ds}{\sqrt{\mathcal{G}(\gamma (s))}}+1\right)}\nonumber \\
     && \nonumber \\
     &=& C+ A \cdot \vert \alpha(X,X)\vert,\nonumber
   \end{eqnarray} since $\psi''\leq 0$. Taking in consideration the bounds \eqref{eqHessGP-2} \& \eqref{eqHessGP-3}, the inequality \eqref{eqHessGP-1} yields, ($\phi'(s)=h(s)$),
   \begin{equation}\left[\frac{A}{k\cdot h(s_k)} + 1\right]\vert \alpha (X,X)\vert \geq (1-\sum_{i=1}^{\ell}c_{i}^{2})\frac{h'}{h}(s_k) - \frac{C}{k \cdot h(s_k)}\label{eqHessGP-4}.
   \end{equation}
   Under the hypotheses of item 1. we have that $(h'/h)(s)\geq C_{b}(s)>\sqrt{\vert b\vert}$ and $h(s)\to \infty$ as $s\to \infty$. Moreover, there exists a subspace $V_{x_k}\subset T_{x_k}M$, ${\rm dim}V_{x_k}\geq 2$,  such that if $X\in V_{x_k}$ then $X=a\cdot \grad \rho_{_N}+\sum_{j=1}^{n-\ell-1}b_{j}\cdot \partial/\partial \theta_{j}.$ Therefore, for $X\in V_{x_k}$, $\vert X\vert=1$, we have for $k\gg 1$.
    \begin{eqnarray}\left[\frac{A}{k\cdot h(s_k)} + 1\right]\vert \alpha (X,X)\vert & \geq &  \frac{h'}{h}(s_k) - \frac{C}{k \cdot h(s_k)} \nonumber \\
           &>& \vert b\vert -\frac{C}{k \cdot h(s_k)} \\
           &>& 0.\nonumber
   \end{eqnarray}The proof follows exactly  the steps of  the proof of Theorem \ref{thmMain-intro} and we obtain that
   $\sup_{_M} K_{M}\geq \vert b\vert + \inf_{_N}K_{N}$ if $b<0$ and $\sup_{_M} K_{M}>\inf_{_N}K_{N}$ if $b=0$.
   \vspace{2mm}

To prove item 2., take an orthonormal basis $X_1,...,X_q,...,X_m\in T_{x_k}M$, $$ X_q=a_q\cdot \grad \rho_{_N}+\sum_{j=1}^{n-\ell-1}b_{jq}\cdot \partial/\partial \theta_{j} + \sum_{i=1}^{\ell}c_{iq}\cdot \partial/\partial \gamma_i$$  with $a_q^2+\sum_{j=1}^{n-\ell-1} b_{jq}^{2}+ \sum_{i=1}^{\ell} c_{iq}^{2}=1.$  Tracing the inequality \eqref{eqHessGP-4} to obtain \begin{eqnarray}\left[\frac{A}{k\cdot h(s_k)} + 1\right]\vert H\vert &\geq & (m-\sum_{q=1}^{m}\sum_{i=1}^{\ell}c_{iq}^{2})\frac{h'}{h}(s_k) - \frac{C}{k \cdot h(s_k)}\nonumber \\ \label{eqHessGP-5}
&\geq &(m-\ell) C_b(s_k)- \frac{C}{k \cdot h(s_k)}\\
&>&0\nonumber \end{eqnarray} for $k\gg 1$. If $b=0$ then $C_b(s)=1/s$ then, coupled  with the estimate  $h(s)\geq s \sqrt{s}$, see \cite{bianchini-mari-rigoli}, we deduce that $\sup_{_M}\vert H\vert >0$. And if $b<0$ then $C_b(s)\geq \sqrt{\vert b\vert }>0$, then letting $k\to \infty$ we have $\sup_{_M} \vert H \vert \geq (m-\ell) \sqrt{\vert b\vert }>0$ if $b<0$.
We can see these curvature estimates as geometric applications  of the following extension of the Pigola, Rigoli, Setti \cite[Thm.1.9]{prs-memoirs}.
\begin{theorem}\label{thm-OYBFL}Let $W$ be a complete Riemannian manifold with an Omori-Yau pair $(\mathcal{G}, \gamma)$ for the Hessian $($Laplacian$)$. If $u \in C^{2}(W)$
satisfies
$
\displaystyle{\lim_{x \rightarrow \infty}\frac{u(x)}{\psi(\gamma(x))}} = 0 ,
$ where $
\displaystyle{\psi(t) = \log\left(\int_{0}^{t}\frac{ds}{\sqrt{\mathcal{G}(s)}}+1\right)},
$
then there exist a sequence ${x_{k}} \in M$, $k \in \mathbb{N} $ such that
\begin{equation} \begin{array}{lllllllllll} \vert \grad_{_W} u\vert(x_{k})& \leq&  \displaystyle\frac{A}{k}, &&
 \Hess_{_W} u(x_{k}) & \leq & \displaystyle\frac{C}{k} &&  ( \triangle_{_W} u(x_{k})&\leq & \displaystyle\frac{B}{k})\end{array}\end{equation}
If  $u^{\ast}=\sup_{M}u<\infty$ then  $u(x_{k})\to u^{\ast}$. The constants $A$, $B$ and $C$ come from the Omori-Yau pair $(\mathcal{G}, \gamma)$, see Definition \ref{def:OY-pair}.

\end{theorem}
This result above should be compared with   \cite[Cor. A1.]{prs-jfa}, due to  Pigola, Rigoli, and Setti where they proved an Omori-Yau for quite general operators, applicable to certain unbounded functions  with  growth to infinity faster than ours. However, we could replace the distance function of their result by an Omori-Yau pair. It would be interesting to understand these facts.

 Assume that the  Omori-Yau pair $(\mathcal{G}, \gamma)$ is for the Hessian. The case of the Laplacian is similar. Fix a point $x_0\in M$ such that $\gamma (x_0)>0$ and  define for each $k\in \mathbb{N}$, $g_{k}\colon M \to \mathbb{R}$  by $g_{k}(x) = u(x) - \displaystyle\frac{1}{k}\psi(\gamma(x)) + 1- u(x_0) - \displaystyle\frac{1}{k}\psi (\gamma(x_0))$. We have that $g_{k}(x_0)=1$ and $g_{k}(x)\leq 0$ for $\rho_{_W}(x)={\rm dist}_{_W}(x_0, x)\gg 1$. Thus there is a point $x_k$ such that $g_{k}$ reaches a maximum. This way we find a sequence $x_k\in M$ such that for all $X\in T_{x_{k}}W$
\begin{eqnarray}\Hess_{_W}u(X,X) &\leq & \frac{1}{k}\Hess_{_W}\psi (\gamma)(X,X)\nonumber \\ &=&\frac{1}{k}\left[ \psi''(\gamma)\langle \grad_{_W}\gamma, X \rangle^{2}+ \psi'(\gamma)\Hess_{_W}\gamma (X,X)\right] \nonumber\\
&\leq & \frac{1}{k}\left[\frac{1}{\sqrt{\mathcal{G}(\gamma)}} \, \displaystyle\frac{1}{\displaystyle\left(\int_{0}^{\gamma} \displaystyle\frac{ds}{\mathcal{G}(s)}\,+\,1\right)}\, C \, \sqrt{\mathcal{G}(\gamma)}\, \displaystyle\left(\int_{0}^{\gamma} \displaystyle\frac{ds}{\mathcal{G}(s)}\,+\,1\right) \right]\vert X\vert^{2}\nonumber \\
&=& \frac{C}{k}\vert X \vert^{2}.\nonumber\end{eqnarray}We used that $\psi''\leq 0$ and $\Hess_{_W}\gamma (X,X)\leq  C \, \sqrt{\mathcal{G}(\gamma)}\, \displaystyle\left(\int_{0}^{\gamma} \displaystyle\frac{ds}{\mathcal{G}(s)}+1\right)$.

\begin{eqnarray}\vert \grad_{_W}u\vert &=& \displaystyle\frac{1}{k}\vert  \grad_{_W}\psi(\gamma)\vert\nonumber \\ &\leq & \displaystyle\frac{1}{k}\left[\frac{1}{\sqrt{\mathcal{G}(\gamma)}} \, \displaystyle\frac{1}{\displaystyle\left(\int_{0}^{\gamma} \displaystyle\frac{ds}{\mathcal{G}(s)}+1\right)}\, A \, \sqrt{\mathcal{G}(\gamma)}\, \displaystyle\left(\int_{0}^{\gamma} \displaystyle\frac{ds}{\mathcal{G}(s)}+1\right)\right]\nonumber\\
&\leq & \displaystyle\frac{A}{k}.\nonumber\end{eqnarray}

\subsection{ Omori-Yau pairs and warped products}Let $(N, g_{_N})$ and $(L, g_{_L})$ be complete Riemannian manifolds of dimension $n-\ell$ and $\ell$ respectively and $\xi\colon L \to \mathbb{R}_{+}$ be a smooth function. Let $\varphi \colon M \to L\times_{\xi}N$ be an isometric immersion into the warped product $L\times_{\xi}N=(L\times N, ds^{2}= g_{_L}+ \xi^2g_{_N})$.  The immersed submanifold $\varphi(M)$ is cylindrically bounded if  $\pi_{N}(\varphi(M))\subset B_{N}(r)$, where $\pi_{N}\colon L\times N \to N$ is the canonical projection in the $N$-factor and $B_{N}(r)$ is a regular geodesic ball of radius $r$ of $N$.
Al\'ias and  Dajczer  in the proof of \cite[Thm.1]{alias-dajczer}, showed that if  $ \varphi $ is  proper in $L\times N$ then   the existence of an Omori-Yau pair for the Hessian in $ L $ induces an Omori-Yau pair for the Laplacian on $M$ provided the mean curvature  $ \vert H \vert $ is bounded. We can prove a slight extension of this result.
\begin{lemma}\label{pullbackomoriyau} Let $ \varphi \colon M \rightarrow L\times_{\xi}N  $ be an isometric immersion, proper in the first entry, where $ L $ carries an Omori-Yau pair $(\mathcal{G},\gamma) $ for the Hessian, $ \xi \in C^{\infty}(L) $ is a positive function satisfying
\begin{eqnarray}\label{hipeta}
\vert \grad \log \xi (y)\vert \leq \ln\left(\int_{0}^{\gamma(y)}\frac{ds}{\sqrt{\mathcal{G}(s)}} + 1\right).
\end{eqnarray}
Letting $\varphi(x)=(y(x), z(x))$ and if
\begin{eqnarray}\label{hipcurvmedia}
\vert H(\varphi(x))\vert \leq \ln\left(\int_{0}^{\gamma (y(x))}\frac{ds}{\sqrt{\mathcal{G}(s)}}+1\right),
\end{eqnarray}then $ M$ has an Omori-Yau pair for the Laplacian. In particular,  $M$ holds the  Omori-Yau maximum principle for the Laplacian.
\end{lemma}
 The idea of the proof is presented in \cite{alias-dajczer} and therefore will try to follow the same notation to simplify the demonstration. Let ($\mathcal{G}, \gamma$) the Omori-Yau pair for the Hessian of $L$.
 Assume w.l.o.g. that $ M $ is non-compact and denote $ \varphi(x) = (y(x), z(x)) $. Define $ \Gamma(y,z) = \gamma(y) $  and define $ \vartheta  (x)= \Gamma\circ \varphi = \gamma(y(x)) $. We will show that $ (\mathcal{G},\vartheta) $ is an Omori-Yau pair for the Laplacian in $ M $.
Indeed, let $ q_k \in M $ a sequence such that $ q_k \rightarrow \infty $ in $ M $ as $ k \rightarrow +\infty $. Since $ \varphi $ is proper in the first entry, we have that $ y(q_k) \rightarrow \infty $ in $ L $. Since $\vartheta(q_k)= \gamma (y(q_k)) $ we have  $ \vartheta(q_k) \rightarrow \infty $ as $ q_{k} \rightarrow \infty $ in $M$.
\vspace{2mm}

We have that
\begin{eqnarray}\label{gradGama}
\grad_{_{L\times_{\xi} N}} \,\Gamma(z,y) = \grad_{_L} \gamma (z) .
\end{eqnarray}
Since $ \xi = \Gamma\circ\varphi $, we obtain at $\varphi (q)$
\begin{eqnarray}
\grad_{_{L\times_{\xi} N}}\, \Gamma& = &(\grad_{_{L\times_{\xi} N}} \Gamma)^{T} + (\grad_{_{L\times_{\xi} N}} \Gamma)^{\perp}\nonumber \\
&=& \grad_{_M} \xi + (\grad_{_{L\times_{\xi} N}}  \Gamma)^{\perp}.\nonumber
\end{eqnarray}
By  hypothesis we have
\begin{eqnarray*}
\vert \grad_{_M}\, \xi \vert(q) &\leq & \vert \grad_{_{L\times_{\xi} N}}\, \Gamma\vert (\varphi(q))
= \vert \grad_{_L}\,\gamma\vert(y(q)) \nonumber \\
&&\\
&\leq & \sqrt{\mathcal{G}(\gamma(y(q)))}\left(\int_{0}^{\gamma(y(q))}\frac{ds}{\sqrt{\mathcal{G}(s)}} + 1\right) \nonumber
\end{eqnarray*}
out of  a compact subset of $ M $.
\vspace{2mm}

Let $T,S\in TL$, $X,Y\in TN$ and  $\nabla^{L\times_{\xi} N}$,  $\nabla^{L}$ and $\nabla^{N}$  be the Levi-Civita connections of the metrics  $ds^{2}=g_{L}+ \xi^2g_{N}$, $g_{L}$ and $g_{N}$ respectively.   It is easy to show that $\nabla^{^{L\times_{\xi} N}}_{S}T= \nabla^{^{L}}_{S}T$ and $ \nabla^{^{L\times_{\xi} N}}_{X}T = \nabla^{^{L\times_{\xi} N}}_{T}X = T(\eta)X $  where $ \eta = \log\xi $. Therefore,
$$\begin{array}{lll}\nabla^{^{L\times_{\xi} N}}_{T}\grad_{_{L\times_{\xi} N}} \Gamma &=&  \nabla^{^{L}}_{_T}\grad_{_L} \gamma \\
&& \\
\nabla^{^{L\times_{\xi} N}}_{X}\grad_{_{L\times_{\xi} N}} \Gamma &=& \grad_{_L} \gamma(\eta)\,X .
\end{array}$$
Hence,  $$\begin{array}{lllll}\Hess_{\!_{L\times_{\xi} N}}\Gamma(T,S)& = & \Hess_{\!_L}\gamma(T,S), \,\,\,   \Hess_{\!_{L\times_{\xi} N}}\Gamma(T,X)& =& 0\\
&&&&\\ \Hess_{\!_{L\times_{\xi} N}}\Gamma(X,Y)&=&\langle \grad_{_L} \eta,\grad_{_L} \gamma\rangle\langle  X,Y\rangle.&&\end{array}$$

For any unit vector $ e \in T_{q}M $, decompose $ e = e^L + e^N $, where $ e^L \in T_{y(q)}L $ and $ e^N \in T_{z(q)}N $. Then we have at $\varphi (q)$
\begin{eqnarray*}
\Hess_{\!_{L\times_{\xi} N}}\Gamma(e,e) = \Hess_{\!_L}\gamma(y(q))(e^L,e^L) + \langle \grad_{_L} \gamma,\grad_{_L} \eta\rangle(y(q))\vert e^N\vert^2 .
\end{eqnarray*}
On the other hand, $\Hess_{\!_M} \xi(q)(e,e)= \Hess_{\!_{L\times_{\xi} N}}\Gamma(e,e)+ \langle \grad_{L\times_{\xi} N} \Gamma,\alpha(e,e)\rangle$. Therefore,
\begin{eqnarray}\label{expression1}
\Hess_{\!_M} \xi(q)(e,e) &=& \Hess_{\!_L}\gamma(e^L,e^L) + \langle \grad_{_L} \gamma,\grad_{_L} \eta\rangle(z(q))\vert e^P\vert^2  \nonumber \\
&& \\
&&+ \,\,\langle \grad_{_L} \gamma,\alpha(e,e)\rangle.\nonumber
\end{eqnarray}

However,
\begin{eqnarray}\label{expressionI}
\Hess_{\!_L}\gamma  \leq \sqrt{\mathcal{G}(\gamma)}\left(\int_{0}^{\gamma}\frac{ds}{\sqrt{\mathcal{G}(s)}} + 1\right),
\end{eqnarray}
out of  a compact subset of $ L $. By hypothesis, see \eqref{hipeta},
\begin{eqnarray}\label{expressionII}
\langle \grad_{_L} \gamma,\grad_{_L} \eta\rangle(y(q)) \!&\leq & \!\vert \grad_{_L}\gamma\vert \cdot \vert \grad_{_L}\eta\vert \nonumber \\
&\leq & \!\sqrt{\mathcal{G}(\gamma)}\left(\int_{0}^{\gamma}\frac{ds}{\sqrt{\mathcal{G}(s)}} + 1\right)\ln\left(\int_{0}^{\gamma}\frac{ds}{\sqrt{\mathcal{G}(s)}} + 1\right).
\end{eqnarray}

Considering  $(\ref{expressionI})$, $(\ref{expressionII})$ and  $(\ref{expression1})$ we have that (off a compact set)
\begin{eqnarray*}
\Hess_{\!_M} \xi(q)(e,e) &\leq & C\cdot \sqrt{\mathcal{G}(\gamma)}\left(\int_{0}^{\gamma}\frac{ds}{\sqrt{\mathcal{G}(s)}} + 1\right)\ln\left(\int_{0}^{\gamma}\frac{ds}{\sqrt{\mathcal{G}(s)}} + 1\right)\nonumber \\
&& \\
&& + \langle \grad_{_L} \gamma,\alpha (e,e)\rangle,
\end{eqnarray*} for some positive constant $C$.
Thus, by $(\ref{hipcurvmedia})$ it follows that
\begin{eqnarray*}
\triangle \gamma \leq B\sqrt{G(\gamma)}\left(\int_{0}^{\gamma}\frac{ds}{\sqrt{\mathcal{G}(s)}} + 1\right)\ln\left(\int_{0}^{\gamma}\frac{ds}{\sqrt{\mathcal{G}(s)}} + 1\right)
\end{eqnarray*} for some positive constant $B$.
Concluding that $(\mathcal{G}, \xi)$ is an Omori-Yau pair for the Laplacian in $M$.
The proof of \cite[Thm.1]{alias-dajczer} coupled with Lemma \ref{pullbackomoriyau}   allows us to state the following
technical extension of Alias-Dajczer's Theorem \cite[Thm.1]{alias-dajczer}.
\begin{theorem}[Alias-Dajczer]\label{Thm-AD}Let $ \varphi \colon M \rightarrow L\times_{\xi}N  $ be an isometric immersion, proper in the first entry, where $ L $ carries an Omori-Yau pair $(\mathcal{G},\gamma) $ for the Hessian, $ \xi \in C^{\infty}(L) $ is a positive function satisfying
\begin{eqnarray}\label{hipeta-B}
\vert \grad \log \xi (y)\vert \leq \ln\left(\int_{0}^{\gamma(y)}\frac{ds}{\sqrt{\mathcal{G}(s)}} + 1\right).
\end{eqnarray}
Letting $\varphi(x)=(y(x), z(x))$ and if
\begin{eqnarray}\label{hipcurvmedia-B}
\vert H(\varphi(x))\vert \leq \ln\left(\int_{0}^{\gamma (y(x))}\frac{ds}{\sqrt{\mathcal{G}(s)}}+1\right).
\end{eqnarray}Suppose that $\varphi (M)\subset \{ (y,z): y\in L,\,z\in B_{_N}(r)\}$ then  $$\sup_{M}\xi\vert H\vert \geq (m-\ell)C_{b}(r), $$ where $b=\sup_{B_{_N}(r)}K^{\rad}_{_N}$.
\end{theorem}
\begin{remark}The Theorems \ref{thmMain-intro} \& \ref{thmMain-2} should have  versions for $\phi$-bounded submanifold of  warped product $L\times_{\xi}N$. Specially interesting should be the Jorge-Koutrofiotis Theorem in this setting. We leave to the interested reader to pursue it.
\end{remark}

As a last application of Theorem \ref{thm-OYBFL}, let $ N^{n+1} = I\times_{\xi}P^n $ the product manifold endowed with the  warped product metric, $ I \subset \mathbb{R} $ is a open interval, $ P^n $ is a complete Riemannian manifold and $ \xi \colon I \to \mathbb{R}_{+} $ is a smooth function. Given an isometrically immersed hypersurface  $ \varphi \colon M^{n} \to N^{n+1} $, define $ h : M^{n} \to I $ the $ C^{\infty}(M^{n}) $ height function by setting $ h = \pi_{I}\circ \varphi $, where $\pi_{I}\colon I\times P\to I$ is a projection. This result below is a technical extension of   \cite[Thm.7]{alias-impera-rigoli} its proof is exactly as there, we just relaxed the hypothesis guaranteeing  an Omori-Yau sequence.

\begin{theorem}\label{hypersufacetheorem} Let $ \varphi \colon M^{n} \to N^{n+1} $ be an isometrically immersed hypersurface. If  $ M^n $ has an Omori-Yau pair $(\mathcal{G}, \gamma)$ for the Laplacian and the height function $ h $ satisfies
$
\displaystyle{\lim_{x \rightarrow \infty}\frac{h(x)}{\psi(\gamma(x))}} = 0$ then
\begin{equation}
\sup_{M^n}\vert H \vert \geq \inf_{M^n}\mathcal{H}(h) ,
\end{equation}
with $ H $ being the mean curvature and $ \displaystyle{\mathcal{H}(t) = \frac{\rho'(t)}{\rho(t)}} $.
\end{theorem}

\vspace{2mm}

\noindent \textbf{Acknowledgements.} We want to express our gratitude to our friend  Newton Santos for their suggestions along the preparation on this paper.



\end{document}